\documentclass [twoside]{article}
\input amssym.def
\input amssym.tex

\newtheorem{thm}{Theorem}[section]
\newtheorem{prop}[thm]{Proposition}
\newtheorem{cor}[thm]{Corollary}
\newtheorem{lem}[thm]{Lemma}
\newtheorem{remark}{Remark} 
\newtheorem{dfn}{Definition} 

\def\bi{{\bf i}}
\def\bI{{\bf J}}

\def\bz{{\bf z}}
\def\b1{{\bf 1}}

\def\bT{{\bf T}}

\def\bR{{\bf R}}

\def\bZ{{\bf Z}}
\def\bC{{\bf C}}

\def\Blat{\mbox{\it \raise2pt\hbox{"}\kern-2pt H}}
\def\lvup{\rlap{\ ${}^{q\atop{\hbox{${}^{\vee}$}}}$}\cdots}
\begin{document}
\begin{center}
{\center{\Large{\bf
Hodge structure of fibre integrals associated to  the affine hypersurface
in a torus} }}

 \vspace{1pc}
{ \center{\large{ Susumu TANAB\'E }}}
\end{center}

\noindent
\begin{center}
 \begin{minipage}[t]{10.2cm}
{\sc Abstract.} {\em
We calculate the fibre integrals of the affine hypersurface
in a torus in the form of their Mellin transforms.
Especially, our method works efficiently for an affine hypersurface
defined by a so called ``simpliciable'' polynomial.
The relations between poles of Mellin transforms of fibre integrals,
the mixed Hodge
structure of the cohomology of the
hypersruface, the hypergeometric differential equation,
and the Euler characteristic of fibres are clarified.
}
 \end{minipage} \hfill
\end{center}
{
\center{\section{Introduction}}
}

In this note we propose a simple method to calculate concretely
fibre integrals asociated to the affine hypersurface in a torus.
We establish an expression of the position of poles of the Mellin
transform with the aid of the mixed Hodge structure of an
hypersurface $Z_f$ defined by a $\Delta-$regular polynomial
explained by  V. Batyrev \cite{Baty1}. The trial to relate the
asymptotic behaviour of a fibre integral with the Hodge structure
of the fibre variety goes back to \cite{Var} where Varchenko
established the equivalence of the asymptotic Hodge structure and
the mixed Hodge structure in the sense of Deligne-Steenbrink for
the case of plane curves  and (semi-)quasihomogneous
singularities.
%



%
\par Later on, several authors (\cite{LS},
\cite{Sab1},\cite{Sab2},\cite{Sab3}) have pursued studies on the
the asymptotic behaviour of fibre integrals in making use of the
Mellin transforms. Their main idea consists in the fact, that
it is possible to visualize the asymptotic behaviour (i.e. the filtration)
of fibre integrals by means of
the poles of Mellin transform. Especially in the case of
complete intersection singularities, the advantage of this method
is quite clear. Let us remark also that not only poles of the Mellin transform
but its zeros play role in the calculus of the global monodromy of the fibre integrals
(e.g. see Proposition ~\ref{prop53}).
The relation between the poles of the Mellin transform and the mixed
Hodge structure has been explained for  examples of
isolated complete intersections of space curve type in \cite{Tan}.
 \par In this note, we illustrate the clarity of this approach
in taking the example of a hypersurface in a torus defined by
so called simpliciable polynomial
(see Definition ~\ref{dfn1}). It  serves as an introduction to the author's
main paper in preparation \cite{Tan2} where he establishes
similar results
for the fibre integrals of  the
``simpliciable''
complete intersection singularities.
\par
I would like to express my gratitude to Prof. M.Oka and the staffs of Tokyo
Metropolitan University who gave me an occasion to discuss on this subject
and to formulate this hypersurface version of the calculus.

{
\center{\section{Hodge structure of the cohomology group of a hypersurface
in a torus}}
}
In this section we review fundamental notions on
the Hodge structure of the cohomology group of a hypersurface
in a torus after \cite{Baty1}, \cite{DX1}.

 \vspace{1.5pc}
\footnoterule

\footnotesize{AMS Subject Classification: 14M109primary), 32S25, 32S40
(secondary).

Key words and phrases: affine hypersurface,
Hodge structure, hypergeometric function.}

 \footnotesize{partially supported by Max Planck
Institut f\"ur Mathematik.}
\normalsize

\par Let $\Delta$ be a convex $n-$dimensional convex polyhedron in
${\bR}^n$ with all vertices in $\bZ^n.$
Let us define a ring $S_\Delta \subset$
$\bC[x_1^{\pm}, \cdots, x_n^{\pm}]$ of
the Laurent polynomial ring as follows:
$$S_\Delta:= \bC \cup \bigcup_{\frac{\vec \alpha}{k} \in \Delta,
\exists k \geq 1} \bC \cdot x^{\vec \alpha}.\leqno(1.1) $$

\par
We denote by $\Delta(f)$ the convex hull of the set ${\vec \alpha}
\in supp(f)$ and call it the Newton polyhedron of a Laurent
polynomial  $f(x).$ We introduce the following Jacobi ideal:
$$ J_{f,\Delta}= \bigl<x_1 \frac{\partial f}{\partial x_1},
\cdots,  x_n \frac{\partial f}{\partial x_n}\bigr>\cdot S_{\Delta(f)}.
\leqno(1.2)$$
Let $\tau$ be a $\ell-$dimensional face of $\Delta(f)$
and define
$$ f^\tau(x)= \sum_{\vec \alpha \in \tau \cap supp(f)}a_{\vec \alpha}
x^{\vec \alpha},
\leqno(1.3)$$
where $f(x) =\sum_{\vec \alpha \in supp(f)}a_{\vec \alpha}
x^{\vec \alpha}.$
The Laurent polynomial $f(x)$ is called $\Delta$ regular, if
$\Delta(f)=\Delta$ and for every $\ell-$dimensional face $\tau
 \subset \Delta(f)$ ($\ell >0$) the polynomial equations:
$$ f^\tau(x)= x_1 \frac{\partial f^\tau}{\partial x_1}=
\cdots=  x_n \frac{\partial f^\tau}{\partial x_n}=0,$$
have no common solutions in $\bT^n = (\bC^\times)^n.$

\begin{prop}
Let $f$ be a Laurent polynomial such that $\Delta(f)=\Delta.$
Then the following conditions are equivalent.
\par
(i) The elements $x_1 \frac{\partial f}{\partial x_1},
\cdots,  x_n \frac{\partial f}{\partial x_n}$ gives rise to a regular sequence
in $S_{\Delta(f)}$
\par
(ii)
$$ dim \bigl(\frac{ S_\Delta}{J_{f,\Delta}}\bigr)= n! vol(\Delta).$$
\par
(iii) $f$ is $\Delta-$regular.
\end{prop}

It is possible to  introduce a filtration on $S_\Delta,$
namely $\vec \alpha \in S_k $ if and only if
$ \frac{\vec \alpha}{k} \in \Delta.$
Consequently we have an increasing filtration;
$$\bC \cong \{0\}=S_0 \subset S_1 \subset \cdots \subset S_n \subset \cdots,$$
that induces a decreasing filtration on  $\frac{ S_\Delta}{J_{f,\Delta}}:$
$$ F^n\bigl(\frac{ S_\Delta}{J_{f,\Delta}}\bigr) \subset
F^{n-1}\bigl(\frac{ S_\Delta}{J_{f,\Delta}}\bigr) \subset \cdots \subset
F^0 \bigl(\frac{ S_\Delta}{J_{f,\Delta}}\bigr).$$
This is called the Hodge filtation of  $\frac{ S_\Delta}{J_{f,\Delta}}.$
It is worthy to remark here that the Hodge filtration ends up with $n-$th
term.
\par
Let us remind us of the notion of Ehrhart polynomial:
\begin{dfn}
Let $\Delta$ be an $n-$dimensional convex polytope. Denote the
Poincar\'e series of graded algebra $S_\Delta$ by
$$P_\Delta(t)= \sum_{k \geq0} \ell(k\Delta)t^k,$$
$$Q_\Delta(t)= \sum_{k \geq0} \ell^\ast(k\Delta)t^k,$$
where $\ell(k\Delta)$ (resp.$\ell^\ast(k\Delta)$ )
represents the number of integer points
in $k\Delta.$ (resp. interior integer points in $k\Delta.$ )
Then
$$\Psi_\Delta(t)= \sum_{k = 0}^n \psi_k(\Delta)t^k= (1-t)^{n+1}P_\Delta(t),$$
$$\Phi_\Delta(t)= \sum_{k = 0}^n \varphi_k(\Delta)t^k = (1-t)^{n+1}Q_\Delta(t),$$
are called Ehrhart polynomials which satisfy
$$t^{n+1}\Psi_\Delta(t^{-1})= \Phi_\Delta(t).$$
\end{dfn}
Further, the main object of our study will be the cohomology group
of the hypersurface $Z_f:=\{x \in \bT^n ; f(x)=0\}.$
We have an important isomorphism on the Hodge filtration of
$PH^{n-1}(Z_f).$
\begin{thm}(\cite{Baty1})
For the primitive part $PH^{n-1}(Z_f)$ of $H^{n-1}(Z_f),$
the following isomorphism holds;
$$ \frac{F^iPH^{n-1}(Z_f)}{F^{i+1}PH^{n-1}(Z_f)} \cong
Gr_F^{n-i}\bigl(\frac{ S_\Delta}{J_{f,\Delta}}\bigr)
=\frac{F^i\bigl(\frac{ S_\Delta}{J_{f,\Delta}}\bigr)}
{F^{i+1}\bigl(\frac{ S_\Delta}{J_{f,\Delta}}\bigr)}. \leqno(1.4)$$
Furthermore
$$ dim \;Gr_F^{n-i} \bigl(\frac{ S_\Delta}{J_{f,\Delta}}\bigr)=
\sum_{q \geq 0}h^{i,q}(PH^{n-1}(Z_f))=\psi_{n-i}(\Delta),$$
for $ i \leq n-1.$

\end{thm}

As for the weight filtration, we have the following characterization.
We understand the notion of the stratum of the support of the algebra
$\frac{ S_\Delta}{J_{f,\Delta}}$
in identifying a polynomial $x^{\vec\alpha} \in S_\Delta$ with
$\vec\alpha \in \bZ^n.$
We call $(n-j)-$dimensional stratum of $supp(S_\Delta)$
the set of those points $\vec i$ from $k \Delta,$ $k=1,2,\cdots$
such that $\frac{\vec i}{k}$ is located on the $(n-j)-$dimensional
face of $\Delta$ and not on any $(n-j-1)-$dimensional
face $\Delta' \subset \Delta.$

\begin{thm}
The weight filtration on $PH^{n-1}(Z_f)$ is defined as a decreasing
filtration
$$0=W_{n-2}\subset W_{n-1} \subset \cdots \subset W_{2n-2} = PH^{n-1}(Z_f),$$
such that $W_{n+i-1}\cong$$\{$ the integer points located on the
strata with dimension $\geq (n-i)$
of $supp\bigl(\frac{ S_\Delta}{J_{f,\Delta}}\bigr)$
but not on the $(n-i-1)-$ dimensional stratum.$\}$ for
$0 \leq i \leq n-2.$
\end{thm}
This theorem is an easy consequence of the Theorem 8.2
\cite{Baty1}. First we notice that the following
exact sequence takes place,
$$ 0 \rightarrow H^n(\bT) \rightarrow H^n(\bT \setminus
Z_f) \stackrel{Res}{\rightarrow} H^{n-1}(Z_f)
\rightarrow 0.$$
The Poincar\'e residue mapping $Res$ gives a morphism of
mixed Hodge structure of the Hodge type $(-1,-1),$
$$ Res(F^j\;H^n(\bT \setminus Z_f)) = F^{j-1}\;H^{n-1}(Z_f),
\;\;Res(W_j\;H^n(\bT \setminus Z_f)) = W_{j-2}\;H^{n-1}(Z_f).$$
Thus we have,
$$ 0 \rightarrow W_{n+i}\;H^n(\bT) \rightarrow
W_{n+i}H^n(\bT \setminus
Z_f) \stackrel{Res}{\rightarrow} W_{n+i-2}H^{n-1}(Z_f)
\rightarrow 0,$$
for $i=2, \cdots , n-1$ where
$$W_{2n-1}\;H^n(\bT)= \cdots = W_{n-1}\;H^n(\bT)=0,
\leqno(1.5)$$
and $dim \; W_{2n}\;H^n(\bT)=1.$
In view of the equality $(1.5)$ the Poincar\'e residue mapping
$Res$ gives an isomorphism
$$ Res:W_{n+i}H^n(\bT \setminus
Z_f) \stackrel{Res}{\rightarrow} W_{n+i-2}H^{n-1}(Z_f),$$
for $i=1, \cdots ,n-1.$ The algebraic structure of the space
$W_{n+i}H^n(\bT \setminus Z_f),$ $i=1, \cdots ,n-1$
has already been established
by Theorem 8.2 \cite{Baty1}.

Further in the course of this paper we identify the element
$x^{\vec\alpha} \in S_\Delta$ with $\frac{x^{\vec\alpha} dx}{df}$
representing an element of  $H^{n-1}(Z_f).$
{
\center{\section{Preliminary combinatorics}}
}
Let us consider a polynomial
$$f(x)= \sum _{1 \leq i \leq M}x^{\vec\alpha(i)} \leqno(2.1)$$
with $M \geq N+1.$ Here $\vec\alpha(i)$ denotes the multi-index
$$\vec\alpha(i)=(\alpha^i_1, \cdots, \alpha^i_N ) \in {\bf Z}^{N}.$$
In the case when $M>N$ we associate to $f(x)$ another polynomial in
$M-1$ variables $f^\sigma(x,x')$
$$f^\sigma(x,x')=\sum_{i=1}^{M-N-1}x_i'x_1^{\alpha({\sigma(i)})}+
\sum_{j=M-N}^M
x_i^{\alpha({\sigma(j)})}
\leqno(2.2)$$
with $\sigma \in {\mathsf{S}}_M,$ the permutation group of $M$
elements. Here we used the notation of the
multi-index:
$$\vec\alpha(\sigma(i))=(\alpha^{\sigma(i)}_1, \cdots, \alpha^{\sigma(i)}_N )
\in {\bf Z}^{N}.$$
In this situation, the expression
$ u(f^\sigma(x,x')+s)$
is a polynomial depending on
$(M+1)$ variables
$(x_1, \cdots, x_N, x'_{1}, \cdots, x'_{M-N-1},s,u).$
Further we shall assume
$$supp(f^\sigma) \cap int(\Delta(f^\sigma))= \emptyset,\leqno(2.3)$$
for each $\sigma$ under question.
Here $\Delta(f^\sigma)$ denotes
the Newton polyhedron of $f^\sigma(x,x').$
\begin{remark}
A polynomial that depends on $(M+1)-$variables and contains
$(M+1)$ monomials is called of Delsarte type. Jean Delsarte
proposed to study algebraic cycles on the hypersurface defined by
a polynomial of this class.
\end{remark}

Let us introduce new variables $T_1, \cdots T_{M+1}$:
$$ T_1 = u x_1' x_1^{\vec\alpha (\sigma(1))},T_2 = u x_2'x_2^{\vec
\alpha(\sigma(2))}, \cdots \leqno(2.4)$$
$$
T_{M-N-1} = u x_{M-N-2}'x^{\vec\alpha(\sigma(M-N-1))},
T_{M-N} = u x^{\vec\alpha(\sigma(M-N))}, \cdots,  T_{M+1}=us.$$
To express the situation in a compact form, we use the following notations:
$$ \Xi := ^t(x_1, \cdots, x_N,  x_1', \cdots,  x'_{M-N-1},u,s),\leqno(2.5)$$
$$Log\; T := ^t(log\; T_1, \cdots, log\; T_{M+1})= ^t(\tau_1, \cdots,
\tau_{M+1}),\leqno(2.6)$$
$$Log\; \Xi := ^t(\log\; x_1, \cdots, \log\; x_N, \log\;x_1', \cdots,
\log\; x'_{M-N-1},\log\;u,, \log\;s). \leqno(2.7)$$

In making use of these notations , we have the relation
$$ \tau_1 = log\;u + log\; x_1' + <\vec\alpha(\sigma(1)), log\; x >,\cdots,
\leqno(2.8)$$
$$\tau_{M-N-1} = log\; u +
log\; x_{M-N-1}' + <\vec\alpha(\sigma(M-N-1)), log\; x >,$$
$$
\tau_{M-N} = log\; u + <\vec\alpha(\sigma(M-N)),log\; x >,
\cdots,
\tau_{M+1} = log\; u +log\;s. $$
We can rewrite the relation (2.8) with the aid of a matrix
${\sf L^{\sigma}} \in End({\bf Z}^{M+1}),$
as follows:
$$ Log\; T= {\sf L^{\sigma}}\cdot Log\; X. \leqno(2.9)$$
where
$${\sf L^{\sigma}}= \left [\begin {array}{cccccccccc} {\it
\alpha_1^{\sigma(1)}}& \cdots&
\alpha_N^{\sigma(1)} &1&0&0&\cdots &0&0&1
\\\noalign{\medskip}\alpha_1^{\sigma(2)}&\cdots &\alpha_N^{\sigma(2)}&
0&1&0 &\cdots &0&0&1\\
\vdots&\cdots&\vdots
&0&0&1&\cdots&0&0&1\\
\\\noalign{\medskip}\alpha_1^{\sigma(M-N-1)}&\cdots &\alpha_N^{\sigma(M-N-1)}&
0&0&0 &\cdots &1&0&1\\
\\\noalign{\medskip}\alpha_1^{\sigma(M-N)}&\cdots &\alpha_N^{\sigma(M-N)}&
0&0&0 &\cdots &0&0&1\\
\vdots&\cdots&\vdots
&\vdots&\vdots&\vdots&\cdots&\vdots&\vdots&\vdots\\
\\\noalign{\medskip}\alpha_1^{\sigma(M)}&\cdots &\alpha_N^{\sigma(M)}&
0&0&0 &\cdots &0&0&1\\
0&\cdots &0&
0&0&0 &\cdots &0&1&1\\
\end {array}\right ],
\leqno(2.10)$$

Further we shall assume that the determinant of the matrix ${\sf L^{\sigma}}$
is positive. This assumption is always satisfied without loss of generality,
if we permute certain column vectors of the matrix, which evidently corresponds
to the change of positions of variables $x.$ We denote the determinant
by $\gamma^\sigma = det(L^{\sigma}).$ The row vectors of $L^{\sigma}$
will be denoted by $\vec e_1^{\sigma}, \cdots, \vec e_{M+1}^{\sigma}. $
 Later we will make use of the notation of variables
$X:= (X_1 , \cdots X_{M-1}):= (x_1, \cdots, x_N,$$
x_1', \cdots,$
$ x'_{M-N-1})$ and that of the polynomial $f^\sigma(x,x')=f^\sigma(X).$
\begin{dfn}
We call that polynomial $f(x)$ is simpliciable if for every
$\sigma \in {\mathsf S}_M,$ $det(L^{\sigma})=\gamma^\sigma \not =0.$
\label{dfn1}
\end{dfn}
For $\tau \subset
\Delta(f^\sigma)$ we denote by $\Sigma(\tau)$
a $(dim\; \tau +1)-$ dimensional simplex consisting all segments connecting
$\{0\}$ and a point of $\tau.$
Let us define a graded algebra
$$ S_\tau := \bigcup_{\frac{\alpha}{k} \in \Sigma(\tau),
\exists k \geq 1} {\bf C}X^\alpha.
\leqno(2.11)$$
and a polynomial
$$f^{\sigma,\tau}(X):= \sum_{\alpha \in supp(f^\sigma) \cap \tau} X^\alpha
\leqno(2.12).$$

\begin{lem}
If $f(x)$ is a simpliciable polynomial, then $f^\sigma(X)$
is $\Delta (f^\sigma)-$ regular.
\par
{\bf Proof}
The condition $det(L^{\sigma})=\gamma^\sigma \not =0$ yields that
$X_1 \frac{\partial f^{\sigma, \tau}}{\partial X_1}, $
$X_2 \frac{\partial f^{\sigma, \tau}}{\partial X_2}, $ $\cdots$
$X_{M-1} \frac{\partial f^{\sigma, \tau}}{\partial X_{M-1}}, $
form a regular sequence in $S_\tau$ for any face $\tau \subset
\Delta(f^\sigma).$ {\bf Q.E.D. }
\end{lem}

{
\center{\section{Mellin transforms}}
}
In this section we proceed to the calculation of  the Mellin transform of
the fibre integrals associated to the hypersurface
$Z_{f^\sigma+s} =\{ X \in \bT^{M-1} ; f^\sigma(X)+s=0\}$
defined by a simpliciable
polynomial.
First of all we consider the fibre integral taken along the fibre
${ \gamma(s)}$
$\in H_{M-2}(Z_{f^\sigma+s})$  as follows,
$$ I^\sigma_{X^{\bI}, \partial \gamma} (s): =
\int_{ \gamma(s)}\frac
{X^{\bI-\b1}dX }{df^\sigma(X)}=
\frac{1}{2\pi \sqrt -1}
\int_{\partial \gamma(s)}\frac
{X^{\bI}dX }{(f^\sigma(X)+s)X^\b1}\leqno(3.1)$$
where ${\partial \gamma(s)}$
$\in H_{M-1}(\bT^{M-1} \setminus Z_{f^\sigma+s})$ is a cycle
obtained after the application of  $\partial,$ Leray's coboundary operator.
Here $X^\b1 = X_1 \cdots X_{M-1},$ $X^{\bI}= X_1^{i_1} \cdots,
 X_{M-1}^{i_{M-1}}.$ See the book by V.A.Vassiliev on ramified integrals  
for the Leray'scoboundary operator.

The Mellin transform of
$I^\sigma_{X^{\bI}, \partial \gamma} (s)$ is defined by the following integral:
$$M^\sigma_{X^{\bI}} (z):=  \int_{\Pi} (-s)^z I^\sigma_{X^{\bI}, \partial \gamma} (s)
\frac{ds}{s}. \leqno(3.2)$$ Here $\Pi$ stands for a cycle in $\bC$
that avoids the poles of $I_{X^{\bI}, \partial \gamma} (s).$
 We assume that
on the set $\partial \gamma^\Pi:= \cup_{s \in \Pi}
\partial \gamma(s), $ $\Re (f^\sigma(X)+s) >0.$
We denote by ${\mathcal L}_q(\bI,z)$ the inner product of $(\bI,z,1)$
with the $q-$th column vector of 
$(\tt L^\sigma)^{-1}.$ 
Let us deform the integral (3.2) in making use of the definition
(3.1):
$$M^\sigma_{X^{\bI}} (z)=
\int_{ {\bf R}_- \times \partial \gamma^\Pi} e^{u(f^\sigma(X)
+ s)} X^{\bI} u (-s)^{z} \frac{du}{u} \wedge \frac{dX}{X^\b1}
\wedge \frac{ds}{s} \leqno(3.3) $$
$$= \frac{1}{\gamma^\sigma}\int_{(L^\sigma)_\ast
({\bf R}_-  \times \partial \gamma^\Pi)} e^{\Psi(T)}
\prod_{q=1}^M T_q^{{\mathcal L}_q(\bI,z)} \prod_{q=1}^M \bigwedge
\frac{dT_q}{T_q},$$ with
$$\Psi(T) = T_1(X,u) + \cdots + T_M(X,u)+T_{M+1}(s,u)= u(f^\sigma(X) + s) \leqno(3.4)$$
where each term $T_i(X,u)$ $(1 \leq i \leq M)$ represents a
monomial term of variables $X,u$  of the polynomial $(3.4)$ while
$T_{M+1}(s,u)=su.$ By virtue of the simple structure of the matrix
$\tt L^\sigma$ $(2.10)$, we can consider the simplex polyhedron
$\tau_q^\sigma \in \bR^{M-1}$ defined as $\bigl<    \vec
e_1^\sigma , \lvup ,\vec e_{M+1}^\sigma \bigr>,$  $ 1 \leq q \leq
M+1 $ where we identify $\vec e_i^\sigma  \in \bZ^{M+1}$ with that
of $\bZ^{M-1}$ after ignoring the last two entries. It means that
we identify $\vec e_i^\sigma$ with the $i-$th row vector of the
matrix $\tt L^\sigma$ of which one removes the last two columns
$^t (0,0, \cdots,0,1),^t (1,1, \cdots ,1) \in \bZ^{M+1}.$ The
chain $\partial \gamma^{\Pi} \times {\bf R}_- $ can be
deformed in $\bC^M$ so far as it does not encounter the
singularity of the integrand.
\begin{prop}
1) The Mellin transform
$M^\sigma_{X^{\bI}} (z)$ of the fibre integral associated to the
simpliciable polynomial $f^\sigma(X)$ has the following form.
$$M^\sigma_{X^{\bI}} (z)= g(z)  \prod_{q=1}^M \Gamma({\mathcal L}_q(\bI,z)),
1 \leq q \leq M+1, \leqno(3.5)$$ where $g(z)$ is a rational
function in $e^{\frac{\pi i z} {\gamma^\sigma}}$ with
$\gamma^\sigma = (M-1)!vol(\Delta(f^\sigma)).$ The linear function
in $(\bI,z),$
$${\mathcal L}_q(\bI,z) = ^t (\bI, z, 1) \vec w_q^\sigma =
\frac{<\vec \alpha_q^\sigma,\bI> + B_q^\sigma z+
C_q^\sigma}{\gamma^\sigma}, \leqno(3.6)$$ where $\vec w_q^\sigma$
is the $q-$th column vector of the matrix $({\sf
L^{\sigma}})^{-1}$.
\par
2) The $M+1$ linear functions ${\mathcal L}_q(\bI,z)$
are classified into the following  three groups.
$${\mathcal L}_{M+1}(\bI,z) = \frac{B^\sigma_{M+1}}
{\gamma^\sigma}z= \frac{\gamma^\sigma}{\gamma^\sigma}z = z.
 \leqno(3.7)_1$$
For $q $ such that $\vec w_q^\sigma=  B^\sigma_q(\vec v_q^\sigma, 1,-1)$
for some $\vec v_q^\sigma \in \bZ^{M-1}, $and $B^\sigma_q \not =0, $
$${\mathcal L}_q(\bI,z) =
\frac{B_q^\sigma(<\vec v_q^\sigma,\bI> + z-1)}{\gamma^\sigma}.
 \leqno(3.7)_2$$

For $q $ such that $\vec w_q^\sigma=  (\vec v_q^\sigma, 0,0)$
for some $\vec v_q^\sigma \in \bZ^{M-1},$ and $B^\sigma_q =0, $
$${\mathcal L}_q(\bI,z) =
\frac{(<\vec v_q^\sigma,\bI>)}{\gamma^\sigma}.
\leqno(3.7)_3$$
Here the case $(3.7)_3$ corresponds to such $q$ that $dim\;\tau_q^\sigma <M-1.$
\par
3)
$$ |B_q^\sigma|= (M-1)!vol(\tau_q^\sigma). \leqno(3.8)$$
\par
4) For $\bI \in \tau_q^\sigma \cap \Delta (f^\sigma),$ with
$dim\;\tau_q^\sigma= M-1,$
$\tau_q^\sigma \not = \Delta (f^\sigma),$
$$ \bigl< \vec v_q^\sigma, \bI \bigr> = 1.$$
$$ \bigl< \vec v_{M+1}^\sigma, \bI \bigr> = 0.$$
\label{prop31}
\end{prop}
{\bf Proof}
1) The definition of the $\Gamma-$ function sounds as follows;
$$ \int_{\bar {\bf R}_-} e^{T}(-T)^\sigma \frac{dT}{T} = (1- e^{2 \pi i
\sigma})\int_{\bR_-} e^{T}(-T)^\sigma \frac{dT}{T} = (1- e^{2 \pi
i \sigma}) \Gamma(\sigma),$$ for the unique nontrivial  cycle
$\bar {\bf R}_-$ turning around $T=0$ that begins and returns to
$\Re T \rightarrow - \infty.$ We apply it to the integral $(3.3)$
and get $(3.5).$ We consider an action on the chain $C_a=\bar {\bf R}_-$
or ${\bf R}_- $ on the complex $T_a$ plane, 
$\lambda:C_a \rightarrow \lambda(C_a)$
defined by the relation,
$$ \int_{\lambda(C_a)} e^{T_a}T_a^{\sigma_a} \frac{dT_a}{T_a}
= \int_{(C_a)} e^{T_a}(e^{2\pi  \sqrt -1 }T_a)^{\sigma_a} 
\frac{dT_a}{T_a}.$$
By means of this action the chain ${\sf L}_\ast (
{\bR_-} \times
\gamma^\Pi ) $ turns out to be homologous to,
$$\sum_{(j_1^{(\rho)},\cdots,j_{M+1}^{(\rho)} )\in [1,
\gamma^\sigma]^{M+1}}
m_{j_1^{(\rho)},\cdots,j_{M+1}^{(\rho)}} 
\lambda^{j_1^{(\rho)}}(\bR_-)\prod_{a'=2}^{M+1}
\lambda^{j_{a'}^{(\rho)}}(\bar {\bR_-} ),$$ with $m_{j_1^{(\rho)},
\cdots, j_{M+1}^{(\rho)}}  \in \bZ. $ This  explains the presence of 
the factor $g(z)=$ $\sum_{(j_1^{(\rho)},\cdots,j_{M+1}^{(\rho)} )\in
[1,\gamma^\sigma]^{M+1}}$ $m_{j_1^{(\rho)},\cdots,j_{M+1}^{(\rho)}}$  
$ e^{2 \pi \sqrt -1 j_1^{(\rho)}{\mathcal
L}_1({\bI, \bz, \zeta})}$ $ \prod_{a'=2}^{M+1}$ $ e^{2 \pi \sqrt -1
j_{a'}^{(\rho)}{\mathcal L}_{a'}({\bI, \bz, \zeta})} (1- e^{2 \pi
\sqrt -1 {\mathcal L}_{a'}({\bI, \bz, \zeta} ) })$ except for the 
$\Gamma-$ function factors.

The points 2)- 5) are reduced to the linear algebra. For
example 3) can be shown, if one remembers the definition of $M$
minors of the matrix $\tt L^\sigma$ calculated in removing the
$M-$th column. 4) If $\bI \in \tau^\sigma_q,$ the vector $\vec
e^\sigma_i$ is orthogonal to $ (\vec v_{M+1}^\sigma, 1,-1)$  for
$i
 \not = q$ and  $\bigl<\vec e^\sigma_q,
B_q^\sigma( \vec v_{M+1}^\sigma, 1,-1)
\bigr> = \gamma^\sigma.$ The result
on the $M-$th and $(M+1)-$st element is explained by the fact
that $\vec e^\sigma_{M+1}= (0,\cdots, 0,1,1)$
is orthogonal to $ (\vec v_{M+1}^\sigma, 1,-1)$  for $ 1 \leq q \leq M.$

{\bf Q.E.D.}
\par
Let us denote the set of such indices
$q$ with strictly positive (resp. strictly negative)
$ B_q^\sigma$ by $I^+
\subset \{1, \cdots, M+1\},$ (resp. by
$I^-
\subset \{1, \cdots, M+1\}$).
The set of indices $q$ for which $ B_q^\sigma=0$
will be denoted by $I^0.$
With these notations, one can formulate the following,
\begin{cor}

1) The Newton polyhedron admits the following representation,
$\Delta (f^\sigma) =$
$\{\vec i \in \bR^M;$
$\bigl<\vec v_{q}^\sigma, \vec i \bigr> \geq 1$
for $q \in I^+,$
$\bigl<\vec v_{\bar q}^\sigma, \vec i \bigr> \leq 1$
for $\bar q \in I^-,$
 $\bigl<\vec v_{q^0}^\sigma, \vec i \bigr> \geq 0$
for $q^0 \in I^0$ $\}. $
\par
2)
We denote by $\chi(Z_{f^\sigma+1})$ the Euler-Poincar\'e
characteristic of the
hypersurface $Z_{f^\sigma+1}=\{X \in \bT^{M-1};f^\sigma(X)+1 \}$
here under the constant $1$ we understand a generic value for $f^\sigma(X).$
The following equality holds,
$$\sum_{q \in I^+} B_q^\sigma=(M-1)!
vol_{M-1} \bigr(\Delta(f^\sigma(X)+1) \bigl)
= (-1)^{M} \chi(Z_{f^\sigma+1}). \leqno(3.9)$$
\par
3) $\sum_{q =1}^{M+1}B_q^\sigma=0. $ In other words,
$$\sum_{\bar{q}\in I^-} B_{\bar{q}}^\sigma =-(\sum_{q \in I^+} B_q^\sigma). \leqno(3.10)$$
\label{cor32}
\end{cor}
{\bf Proof}
1)
After the definition of vectors $\vec v_{1}^\sigma, \cdots,
 \vec v_{M}^\sigma$ we can argue as follows. If $\vec i
$ does not belong to the hyperplane
$\bigl<    \vec e_1^\sigma , \lvup ,\vec e_{M}^\sigma \bigr>,$
then $\bigl<\vec v_{q}^\sigma, \vec i \bigr>= 1+
\frac{\gamma^\sigma}{B_q^\sigma}.$
In the case when $ q \in I^+$ (resp. $\bar{q}\in I^-$)
$\bigl<\vec v_{q}^\sigma, \vec i \bigr>>1$(resp. $\bigl<\vec
v_{q}^\sigma, \vec i \bigr><1$)
that is equivalent to say that all the points $\vec i$
of the Newton polyhedron
$\Delta (f^\sigma)$
satisfy  $\bigl<\vec v_{q}^\sigma, \vec i \bigr> \geq 1  $
for $q \in I^+$ (resp.$\bigl<\vec v_{q}^\sigma, \vec i \bigr> \leq 1  $
for $q \in I^-$).
If $\vec i
\in \bigl<    \vec e_1^\sigma , \lvup ,\vec e_{M}^\sigma \bigr>,$
then $\bigl<\vec v_{q}^\sigma, \vec i \bigr> = 1.$
For $ q^0 \in I^0,$
$\Delta (f^\sigma) \subset
\{\vec i; \bigl<\vec v_{q^0}^\sigma, \vec i \bigr> \geq 0\}, $
because $\bigl<\vec v_{q^0}^\sigma, \vec i \bigr>=1$
for $\vec i
\not \in \bigl<    \vec e_1^\sigma , {\rlap{\ ${}^{q^0
\atop{\hbox{${}^{\vee}$}}}$}\cdots} ,\vec e_{M}^\sigma \bigr>.$
As all possible cases are exhausted by $I^+,I^-,I^0,$
$|I^+|+|I^-|+|I^0|=M +1.$
This yields the statement.
2) Apply the Theorem by \cite{X},\cite{Oka}
on the Euler characteristic. 3) The $(M+1)-$st column vector of $\tt L^\sigma$
is orthogonal to the $M-$th row vector of ${\tt L^\sigma}^{-1},$
$(B_1^\sigma, \cdots, B_{M+1}^\sigma).$

\begin{cor}
Under the above situation, the Mellin inverse of
$M_{X^{\bI}, \gamma}^\sigma (s)$ with properly chosen
periodic function $g(z)$ with period $\gamma^\sigma:$
$$I_{X^{\bI}, \gamma'}^\sigma (s)
= \int_{\check \Pi} g(z)\frac{
\prod_{a\in I^+}\Gamma
\bigl({\mathcal L}_a (\bI , z )\bigr)}
{\prod_{\bar{a}\in I^-}\Gamma \bigl
(1-{\mathcal L}_{\bar{a}}
(\bI,z)\bigr)} s^{-z} dz, \leqno(3.11)$$
defines a convergent analytic function in $ -\pi <arg\; s <\pi.$
\end{cor}
{\bf Proof}
In applying the Stirling's formula
$$ \Gamma(z+1) \sim  (2\pi z)^{\frac 1 2} z^z e^{-z},\;\; \Re \; z \rightarrow +\infty,$$
to the integrand of $(3.11)$, we take into account the relation
$(3.10).$
Here we remind us of the formula $\Gamma(z)
\Gamma(1-z)= \frac{\pi}{sin\;\pi z}.$
As for the choice of the periodic function $g(z)$ one makes
use of N\"orlund's technique \cite{Nor}.
In this way we can choose such $g(z)$
that the integrand is of exponential decay on  $\check \Pi.$
{\bf Q.E.D.}

\par
{\bf Example}
Let us illustrate  the above procedures
by a simple example.
$$f(x)= x_1^5+ x_1^2x_2 + x_1x_2^2+x_2^4.\leqno(3.12)$$
We have 4 possibilities to add a new variable $x_1'$
so that the polynomial (3.12) becomes a simplicial.
$$f^{\sigma_1} (x,x')= x_1'x_1^5+ x_1^2x_2 + x_1x_2^2+x_2^4.$$
$$f^{\sigma_2} (x,x')= x_1^5+ x_1'x_1^2x_2 + x_1x_2^2+x_2^4.$$
$$f^{\sigma_3} (x,x')= x_1^5+ x_1^2x_2 + x_1'x_1x_2^2+x_2^4.$$
$$f^{\sigma_4} (x,x')= x_1^5+ x_1^2x_2 + x_1x_2^2+ x_1'x_2^4.$$
Let us calculate $\tt L^{\sigma_3}$ and  $({\tt L^{\sigma_3}})^{-1}.$
$${\sf L^{\sigma_3}}= \left
[\begin {array}{ccccc}
5& 0& 0& 0&1\\
2& 1& 0& 0&1\\
1& 2& 1& 0&1\\
0& 4& 0& 0&1\\
0& 0& 0& 1&1\\
\end {array}\right ],$$
$$({\sf L^{\sigma_3}})^{-1}= \frac{1}{7}
\left
[\begin {array}{ccccc}
3& -4& 0& 1&0\\
2& -5& 0& 3&0\\
1& -6& 7& -2&0\\
8& -20& 0& 5&7\\
-8& 20& 0& -5&0\\
\end {array}\right ].$$
We have $$ {\mathcal L}_1(\bI,z)=\frac{3i_1 +2i_2 +i_3 + 8
(z-1)}{7},
 {\mathcal L}_2(\bI,z)=\frac{-4i_1-5i_2 -6i_3 -20(z-1)}{7},
 {\mathcal L}_3(\bI,z)=\frac{7i_3}{7},$$
 $$
 {\mathcal L}_4(\bI,z)=\frac{i_1 +3i_2 -2i_3 + 5
(z-1)}{7}, {\mathcal L}_5(\bI,z)=\frac{7z}{7}.$$ Let us denote by
$\vec e_1 =(5,0,0),$$\vec e_2 =(2,1,0),$$\vec e_3 =(1,2,1),$ $\vec
e_4 =(0,4,0),$$\vec e_5 =(0,0,0).$ Then we have
$$ vol(\tau_5)=3! vol(\vec e_1,\vec e_2,\vec e_3,\vec e_4)=7.$$
Similarly $vol(\tau_4)=5,$ $vol(\tau_3)=0,$  $vol(\tau_2)=20,$
 $vol(\tau_1)= 8.$
Remark $\tau_1 +\tau_3+\tau_4 +\tau_5=\tau_2$
(a subdivision of simplex into three simplices)
which yields $ 7+8+5=20.$
The face not affected (see Definition ~\ref{dfn3} below) by $\sigma_3$ is that
spanned by $\vec e_1, \vec e_2,\vec e_4.$

{
\center{\section{
Hodge structure of the fibre integrals
}}
}

Now we can state the relationship between the Hodge structure of
the $PH^{M-2}(Z_{{f^\sigma}})$ and the poles of the Mellin transform
after suitable period function multiplication
$\frac{1}{g(z)}M_{X^{\bI}, \gamma'}^\sigma (z) = \frac{
\prod_{a\in I^+}\Gamma
\bigl({\mathcal L}_a (\bI , z )\bigr)}
{\prod_{\bar{a}\in I^-}\Gamma \bigl
(1-{\mathcal L}_{\bar{a}}
(\bI,z)\bigr)}.$ Here we remind us of the relation
$ \Gamma(z)\Gamma(1-z)= \frac{\pi}{sin \pi z}.$
We will misuse the expression
``the poles of the Mellin transform''  in meaning those of
$ \frac{
\prod_{a\in I^+}\Gamma
\bigl({\mathcal L}_a (\bI , z )\bigr)}
{\prod_{\bar{a}\in I^-}\Gamma \bigl
(1-{\mathcal L}_{\bar{a}}
(\bI,z)\bigr)}.$
\begin{thm}
1) For $X^{\bI} \in Gr_F^p Gr^w_{M-2}PH^{M-2}(Z_{f^\sigma}),$
$0 \leq p \leq M-1,$
the following properties hold
\par
$\bf a)$ $$0 < \bigl<\vec v_{q^0}^\sigma,\bI \bigr >
< M-1-p \;\;  {for}\; q^0\in I^0 $$
$$M-1-p< \bigl<\vec v_{q}^\sigma,\bI \bigr >
< (M-1-p)(1+ \frac{\gamma^\sigma}{B_q^\sigma}) \;\; q \in I^+ $$
$$(M-1-p)(1+ \frac{\gamma^\sigma}{B_{\bar q}^\sigma})
< \bigl<\vec v_{\bar{q}}^\sigma,\bI \bigr > < M-1-p \;\;
 {for}\; {\bar{q}}\in I^- $$
if not $\bigl<\vec v_q^\sigma,\bI \bigr > =0.$
\par
$\bf b)$ The maximal pole of the Mellin transform satisfies;
$$1-(M-1-p)(1+max_{q \in I^+}\frac{\gamma^\sigma}{B_q^\sigma})
< z < 2+p-M.$$
Here the pole is not necessarily a simple pole.
\par
2) For $X^{\bI} \in Gr_F^p Gr^w_{M-1}PH^{M-2}(Z_{f^\sigma}),$
$0 \leq p \leq M-1,$
the following properties hold
\par
$\bf a)$ There exists unique index $ q\in I^+$ such that:
$$\bigl<\vec v_q^\sigma,\bI \bigr >
= M-1-p \;\; $$
\par
$\bf b)$ The maximal pole of the Mellin transform is the simple pole
$$z = 2+p-M.$$
\par
3) For $X^{\bI} \in Gr_F^p Gr^w_{M-2+r}PH^{M-2}(Z_{f^\sigma}),$
$1 \leq r \leq M-3,$
$0 \leq p \leq M-1,$
the following properties hold.
\par
$\bf a)$ There exist $r$ indices $ q_1, \cdots q_{r} \in I^+$ such that:
$$\bigl<\vec v_{q_1}^\sigma,\bI \bigr > =
\bigl<\vec v_{q_2}^\sigma,\bI \bigr >= \cdots=
\bigl<\vec v_{q_{r}}^\sigma,\bI \bigr >
= M-1-p, \;\; $$
but no such  $r+1$ pair of indices  $ q_1, \cdots, q_{r+1}.$
\par
$\bf b)$ The maximal pole of the Mellin transform satisfies;
$$z = 2+p-M,$$
which is of order $\leq r+1$
i.e. there can be cancellation of poles.
\end{thm}
The defect number $(r+1)-$ $\{$order of poles $\}$
will be described in $\S 5.$

Proof of the theorem can be achieved by a combination of Theorems 1.2, 1.3
and the Proposition 3.1, Corollary ~\ref{cor32}. We remember here that
the $\Gamma(z)$
has simple poles at $z=0,-1, -2, \cdots.$
\par

The above theorem mentions about how the Hodge structure of
 $PH^{M-2}(Z_{f^\sigma})$ influences on the poles of the Mellin transform.
How about the original Hodge structure $PH^{N-1}(Z_f)$ ?
To state this relationship, we need to introduce the following notion.
\begin{dfn}
The face $\tau \in\Delta(f)$ is called
``not affected  by $\sigma$'' $\in {\mathsf S}_{M}$
if $\tau \in \Delta(f^\sigma)$ after the extension of $(i_1, \cdots, i_N)
\in\tau  \subset \bR^{N} $
into $\bR^M$
transforming it into the vector $(\bi,0)=(i_1, \cdots, i_N,0,\cdots,0,0)
\in \bR^M.$
\label{dfn3}
\end{dfn}
The face not affected  by $\sigma$
for the polynomial (2.2) is a face (or its sub-face)
spanned by the vertices
 $$\sum_{j=M-N}^M
x_i^{\alpha({\sigma(j)})}$$ i.e. vertices free of $x_i'.$
\begin{thm}
1) For $x^{\bi} \in Gr_F^p Gr^w_{N-1}PH^{N-1}(Z_{f}),$
$0 \leq p \leq N,$ for which $(\bi,0)$ lies in $supp
\bigl(\frac
{S_{\Delta(f^\sigma)}}{J_{f^\sigma, \Delta(f^\sigma)} }\bigr) $
not affected by $\sigma,$
the following properties hold
\par
$\bf a)$ $$0 < \bigl<{\vec v_{q^0}^\sigma},(\bi,0) \bigr >
< N-p \;\;  {for}\; q^0 \in I^0, $$
$$N-p < \bigl<{\vec v_q^\sigma},(\bi,0) \bigr >
< (N-p)(1+ \frac{\gamma^\sigma}{B_q^\sigma} ) \;\;  {for}\; q\in I^+, $$
$$(N-p)(1+ \frac{\gamma^\sigma}{B_{\bar{q}}^\sigma} )
< \bigl<{\vec v_{\bar{q}}^\sigma},(\bi,0) \bigr > < N-p \;\;
 {for}\; {\bar{q}}\in I^-, $$
if not $\bigl <{\vec v_q^\sigma},(\bi,0) \bigr > =0,$
or $\bigl <{\vec v_{\bar{q}}^\sigma},(\bi,0) \bigr > =0.$
\par
$\bf b)$ The maximal pole of the Mellin transform satisfies;
$$ 1-(N-p)(1+ max_{q \in I^+}\frac{\gamma^\sigma}{B_q^\sigma})
< z < 1-N+p.$$
Here the pole is not necessarily a simple pole.
\end{thm}
The proof is straightforward if one applies Theorem 4.1
to $\Delta(f).$ We consider the $N-$dimensional
face  $\tau^\sigma_q \subset \bZ^N $ that is a $N-$dimensional simplex
contained in $\Delta(f).$ One can verify that there exist
$(\bi,0)\in supp
\bigl(\frac {S_{\Delta(f^\sigma)}}{J_{f^\sigma, \Delta(f^\sigma)}}\bigr)$
such that $x^{\bi} \in Gr_F^p Gr^w_{N-1}PH^{N-1}(Z_{f}),$
$0 \leq p \leq N-1$ for the cases $N=2,3,4$ by means of polyhedra
realizing  the formulae 5.11, \cite{DX1}.

\par
We remark the following simple combinatorial fact.
\begin{prop}
For every $x^{\bi} \in Gr_F^p Gr^{w}_{N-1}PH^{N-1}(Z_{f}),$
there exists an element $\sigma \in {\mathsf S}_M$ such that
$x^{\bi}$ is not affected by $\sigma.$
That is to say there    exists $\sigma \in {\mathsf S}_M$ such that
$x^{\bi}  \in S_{\Delta(f)}\cap S_{\Delta(f^\sigma)}.$
\end{prop}
{
\center{\section{
Hypergeometric group associated to the fibre integrals
}}
}

Let us introduce two differential operators of order
$\Delta^\sigma := (M-1)!vol_{M-1}(\Delta(f^\sigma(X)+1))=
|\chi(Z_{f^\sigma +1})|$
$=|I^+|=|I^-|;$
$$P_{\bI}^\sigma(\vartheta_s)=
 \prod_{q\in I^+} \prod_{j=0}^{B_q^\sigma-1}
\bigl({\mathcal L}_q(\bI,-\vartheta_{s})+j \bigr)
\leqno(5.1)$$
$$Q_{\bI}^\sigma( \vartheta_s)
= \prod_{\bar{q}\in I^-}\prod_{j=0}^{-B_{\bar q}^\sigma-1}
\bigl(-{\mathcal L}_{\bar{q}}(\bI,-\vartheta_{s})-j \bigr) , \leqno(5.2)$$
where $I^+, I^- $are those sets of indices introduced in \S 3. We 
have the following theorem as a corollary to the 
Proposition ~\ref{prop31}.

\begin{thm}
The fibre integral
$I_{X^{\bI}, \gamma}^\sigma (s)$ is annihilated by the operator
$$R_{\bI}^\sigma(\vartheta_s)= P_{\bI}^\sigma(\vartheta_s) -
s^{\gamma^\sigma}Q_{\bI}^\sigma(\vartheta_s), \leqno(5.3)_1$$
that is to say
$$[P_{\bI}^\sigma(\vartheta_s) -
s^{\gamma^\sigma}Q_{\bI}^\sigma( \vartheta_s)]I_{X^{\bI}, \gamma}^\sigma 
(s)=0. \leqno(5.4)$$
\end{thm}
It is worthy to remark that  the operator $R_{\bI}^\sigma(\vartheta_s)$
is a push-forward of the Pochhammer hypergeometric operator of order
$\Delta^\sigma,$
$$P_{\bI}^\sigma(\gamma^\sigma \vartheta_t) - t Q_{\bI}^\sigma
(\gamma^\sigma \vartheta_t), \leqno(5.3)_2$$
by the Kummer covering  $t=s^{\gamma^\sigma}.$
In certain cases, the operator $(5.3)_\ast$
turns out to be reducible. Let us introduce the following
set of rational numbers.
$$ C^+(\bI)= \bigcup_{q \in I^+}\bigcup_{0 \leq j \leq B_q^\sigma
-1}  \{ \frac{j} {B_q^\sigma}-
\frac{(<\vec v_q^\sigma,\bI> -1)}{\gamma^\sigma}\}.$$
$$ C^-(\bI)= \bigcup_{{\bar q}
\in I^-}\bigcup_{1 \leq j \leq -B_{\bar q}^\sigma}  
\{ \frac{j} {B_{\bar q}^\sigma}-
\frac{(<\vec v_{\bar q}^\sigma,\bI> -1)}{\gamma^\sigma}\}.$$
$$ C^0(\bI)= C^+(\bI) \cap C^-(\bI).$$
We define a positive integer $\bar \Delta^\sigma =
\sharp |C^+(\bI) \setminus C^0(\bI)|$ $=
\sharp |C^-(\bI) \setminus C^0(\bI)|.$
Then ``the irreducible part''
of  $(5.3)_2$ (i.e. after the division by
operators with rational function solution of type $t^{\alpha^0},$
$\alpha^0 \in C^0(\bI)$)
can be defined as
$$ \bar R_{\bI}^\sigma(\vartheta_t)=
\prod_{\alpha^+ \in C^+(\bI)\setminus C^0(\bI) }(\vartheta_t+
\alpha^+) - t \prod_{\alpha^- \in C^-(\bI)\setminus
C^0(\bI)}(\vartheta_t+ \alpha^-+1),$$ as an operator of order
$\bar \Delta^\sigma$ up to multiplication by a constant to the
variable $"t".$
\par We consider solutions $u_{\ell,m}(t), 1 \leq \ell
\leq \bar \Delta^\sigma,$
to the equation
$$ \bar R_{\bI}^\sigma(\vartheta_t)u_{\ell,m}(t)=0,\leqno(5.5)$$
with the asymptotic behaviour
$$ u_{\ell,m}(t) \cong t^{\rho^\ell_\bI } \;
\sum_{\nu=0}^{m}(log\;t)^\nu A_{\ell, \nu}(t).\leqno(5.5)_1$$
Here
$0 \leq m \leq m_\ell,$
$\sum_{\ell}(m_\ell +1) = \bar\Delta^\sigma,$
$A_\ell(t)$ holomorphic in the neighbourhood
of $t=0.$ Similarly, we consider the asymptotic behaviour
at $t=\infty$ of the solutions to $(5.5)$
$$ v_{\ell,k}(t) \cong (\frac{1}{t})^{\bar \rho^\ell_\bI } \;
\sum_{\mu =0}^{k} (log\;t)^\mu B_\ell(\frac{1}{t}).$$
Here $0 \leq k \leq k_\ell,$  $\sum_{\ell}(k_\ell +1) = \bar \Delta^\sigma,$
$B_\ell(\frac{1}{t})$ holomorphic in the
neighbourhood of $\frac{1}{t}=0.$
Here $m_\ell +1$ (resp.$k_\ell +1$) denotes the multiplicity of
 $-\rho^\ell_\bI$ (resp. $ -\bar\rho^\ell_\bI$)
in the set $C^+(\bI) \setminus C^0(\bI)$ (resp. $C^-(\bI)
\setminus C^0(\bI)$).

Under this situation, we define  characteristic polynomials
of the exponents of solutions to $(5.5)$
at $t=0$
$$X_{0,\bI}({\mathsf t})=\prod_{\ell=1}^{\bar \Delta^\sigma}
({\mathsf t}-e^{2\pi \rho^\ell_\bI \sqrt -1})=
\prod_{\alpha^+ \in C^+(\bI) \setminus C^0(\bI)}
({\mathsf t}-e^{-2\pi \sqrt -1 \alpha^+}),$$
and $t=\infty$
$$X_{\infty,\bI}({\mathsf t})=\prod_{\ell=1}^{\bar \Delta^\sigma}
({\mathsf t}-e^{2\pi \bar \rho^\ell_\bI \sqrt -1}) =
\prod_{\alpha^+ \in C^-(\bI) \setminus C^0(\bI)}
({\mathsf t}-e^{-2\pi \sqrt -1 \alpha^-}).$$
Especially in the case $C^0= \emptyset,$ we have the following
simple formulae.
\begin{cor}
The characteristic polynomials defined above
can be calculated in the following way.
$$X_{0,\bI}({\mathsf t})=\prod_{q \in I^+}
({\mathsf t}^{B^\sigma_q}-e^{-2\pi (1- \bigl
<\vec v_{q}^{\sigma},\bI \bigr >)
\frac{B^\sigma_{q}}{\gamma^\sigma }\sqrt -1}),
\leqno(5.6)_1$$
$$X_{\infty,\bI}({\mathsf t})=\prod_{ q \in I^-}
({\mathsf t}^{-B^\sigma_{\bar q}}- e^{-2\pi
(1- \bigl <\vec v_{\bar q}^{\sigma},\bI \bigr >)
\frac{B^\sigma_{\bar q}}{\gamma^\sigma }\sqrt -1}). \leqno(5.6)_2$$
\end{cor}
For the polynomials introduced in $(5.6)_1,(5.6)_2,$
we introduce two vectors $(A_1, A_2, \cdots, A_{\bar \Delta^\sigma}),$
$ (B_1, B_2,$ $\cdots,$ $ B_{\bar \Delta^\sigma})$
$\in \bC^{\bar \Delta^\sigma}, $
after the following relation:
$$X_{0,\bI}({\mathsf t})={\mathsf t}^{\bar \Delta^\sigma}+
A_1{\mathsf t}^{\bar \Delta^\sigma-1}+
A_2{\mathsf t}^{\bar \Delta^\sigma-2}+ \cdots + A_{\bar \Delta^\sigma},$$
$$X_{\infty,\bI}({\mathsf t})={\mathsf t}^{\bar \Delta^\sigma}+
B_1{\mathsf t}^{\bar \Delta^\sigma-1}+
B_2{\mathsf t}^{\bar \Delta^\sigma-2}+ \cdots + B_{\bar \Delta^\sigma}.$$
Let us denote by
$\omega^i, i=0,1,2, \cdots,\gamma^\sigma -1$
the non-zero singular points of the equation $(5.4)$ i.e.
$\{s \in \bC ;{\prod_{q \in I^+}B^\sigma_q} -
\bigl({\prod_{\bar q \in I^-}B^\sigma_{\bar q}}\bigr)s^{\gamma^\sigma}=0\}.$
\begin{prop}
A representation of the hypergeometric group (global monodromy group)
of the solutions to $(5.5)$ is given by
$$ M_0= h_0^{\gamma^\sigma},
M_{\omega^0} = h_1= (h_0 h_\infty)^{-1}, M_\infty =
h_\infty^{\gamma^\sigma},
M_{\omega^i} = h_\infty^{-i}h_1 h_\infty^i (i=1,2, \cdots,\gamma^\sigma -1),
\leqno(5.7)
$$
for the matrices
$$ h_0=
\left(
\begin{array}{llccll}
0 & 0 & \cdots &0 & -A_{\bar \Delta^\sigma} \\
1 & 0 &  \ddots &0 &-A_{\bar \Delta^\sigma-1} \\
0 &1 & \ddots  &0&-A_{\bar \Delta^\sigma-2} \\
\vdots &\ddots  & \ddots &\vdots &\vdots\\
0 & 0 & \cdots &  1 &-A_1 \\
\end{array} \right),
\leqno(5.8)$$
$$ (h_\infty)^{-1}=
\left(
\begin{array}{llccll}
0 & 0 & \cdots &0 & -B_{\bar \Delta^\sigma} \\
1 & 0 &  \ddots &0 &-B_{\bar \Delta^\sigma-1} \\
0 &1 & \ddots  &0&-B_{\bar \Delta^\sigma-2} \\
\vdots &\ddots  & \ddots &\vdots &\vdots\\
0 & 0 & \cdots &  1 &-B_1 \\
\end{array} \right).
$$
where $M_{\omega^i}$ denotes the monodromy action around the point
$\omega^i \in {\bf CP}^1_s.$
\label{prop53}
\end{prop}
{\bf proof} The monodromies of the solutions annihilated by
$\bar R_{\bI}^\sigma(\vartheta_t)$
are given by $h_0,$ (resp. $h_1, h_\infty$) after
\cite{Lev}.
at $t=0,$ (resp.$t=1,\infty$).
Let us think of a $\gamma^\sigma-$leaf covering
$\tilde {\bf CP}^1_t$ of ${\bf CP}^1_s$
that corresponds to the Kummer covering $s^{\gamma^\sigma} =t.$
In lifting up the path around $t=1$ the first leaf of
$\tilde {\bf CP}^1_s,$ the monodromy $h_1$ is sent
to the conjugation with a path around $t=\infty.$ That
is to say we have $M_{\omega^1}= h_\infty^{-1}h_1 h_\infty.$
For other leaves the argument is similar. {\bf Q.E.D.}

\par
In combining the above result with that of Theorem 4.1, 3),
we get the following.
\begin{cor}
For $X^\bI \in Gr_F^pGr_{M-2+r}^w PH^{M-2}(Z_{f^\sigma}),$
$1 \leq r \leq M-2,$ $0\leq p \leq M-1$
the size of a Jordan cell of the monodromies $M_0$ with unit eigenvalue
arising from the term of the form $(5.5)_1$ with $\alpha^+=
\rho^\ell_\bI$
is $r+1 -  \sharp \{\alpha^+ \in C^0(\bI)
; \alpha^+ \in \bZ\}.$
\end{cor}
{\bf proof}
It is enough to remember the
following relation for a cycle $C$ avoiding $z+\alpha=0$:
$$ (r+1)! \int_C \frac{s^{-z}}{(z+\alpha)^{r+1}} dz=
\int_C s^{-z}[(\frac{d}{dz})^r\frac{1}{(z+\alpha)}] dz$$
$$=\int_C \frac{1}{(z+\alpha)}[(-\frac{d}{dz})^r s^{-z}] dz
=  \int_C \frac{1}{(z+\alpha)} s^{-z} (log\;s)^r dz
= 2 \pi \sqrt -1 s^{\alpha} (log\;s)^r.$$
If the set $C^0(\bI)$ is empty, the order of the poles of
the Mellin transform for  $X^\bI$ $\in Gr_F^pGr_{M-2+r}^w$
$PH^{M-2}(Z_{f^\sigma})$ is $r+1$ after Theorem $4.1, 3) a).$
If $C^0(\bI)$ is not empty, the order of poles is reduced by
$\sharp \{\alpha \in C^0(\bI); \alpha \in \bZ \}.$
{\bf Q.E.D.}

\vspace{\fill}

%

\noindent

\begin{flushleft}
\begin{minipage}[t]{6.2cm}
  \begin{center}
{\footnotesize Indepent University of Moscow\\
Bol'shoj Vlasievskij pereulok 11,\\
 Moscow, 121002,\\
Russia\\
{\it E-mails}:  tanabe@mccme.ru, tanabe@mpim-bonn.mpg.de \\}
\end{center}
\end{minipage}
\end{flushleft}

\end{document}